\documentclass[preprint,review,12pt]{elsarticle}

\usepackage{amssymb}

\journal{journal}
\begin{document}
\begin{frontmatter}
\title{On optimization and solution of roots of a function using Taylor's expansion and fractional derivatives}
\author{Ali Dorostkar$^{1}$ \corref{cor1} \cortext[cor1] {Corresponding author}\\ Ahmad Sabihi$^{2}$}
\ead{Corresponding author: sabihi2000@yahoo.com }
\address{Fractal Group,Isfahan, Iran$^{1}$\\Professor of Mathematics at some of Iranian Universities, Isfahan, Iran $^{2}$}
\begin{abstract}
A method is given for finding roots of a one-variable function using Taylor's expansion of that function and fractional derivative calculated at a suitable tangent point without using Newton's method, but is regarded as a variant of Halley and Newton's one. Several examples regarding polynomials are stated as well. Then, the given method is generalized to functions of several variables belonging to an $n$-dimensional space and one example is given for optimization and solution of a nonlinear system of equations by both our method and Gradient Descent one. A comparison of our method is made with Gradient one for a system of the functions of three variables.  Our given method seems to be much more rapidly than the Newton's one since by finding a suitable point on the function's curve, the number of iterations is to be much less than Newton's iterative steps. We also find order of fractional derivative, which corresponds to equation's found root and compare tangent lines drawn at the root by both fractional and classical derivatives. The methods given in this paper can be used for optimization of function via fractional derivatives of order $\beta$.    
\end{abstract}
\end{frontmatter}
\section{Introduction}
As is well-known, mathematics of convex optimization and finding roots of a function have been investigated more than a century and very interesting subjects in this topic have been studied by many researchers in related recent developments. Some of developments in this area are of more applicable in practice than before. Automatic control systems, signal processing, communications and networks, data analysis, electronic circuit design and modeling, finance and statistics, and etc. are of most important applications of optimization since 1990s. Although, we do not wish to deal with the convex or non-convex optimization methods in this paper, since they have advanced much more than our conception. In this paper and section 2, we present a method similar to that of Newton's one to find roots of an arbitrary equation and then generalize it to functions of several variables. We really make use of the Halley's method given in 1694 (see \cite{HB} and \cite{DFB}).
Several examples are given for clarifying the new method. While presenting the new method of solving an equation, we state another method for finding order of a fractional derivative substituting classical one. The idea of using fractional derivatives in root-finding algorithms is back to Ak\"{u}l et al. \cite{ACT} and Candelario et al.\cite{CCT}, but first of all, we find roots and then make use of fractional derivatives for making a fractional tangent line against classical one and find order of the fractional derivative, which corresponds to equation's root. This is a distinction between our method and other ones.\\  
In section3, we present conclusions of this paper. We suggest to readers some documents to further studies such as: \cite{LRS}, \cite{HBJ}, \cite{HBJC},\cite{KV}, \cite{BA}
\section{Theory}
\subsection{\textbf{The functions of one variable}}
As is well-know, using the Newton's method (\cite{DBP},\cite{GO},\cite{SG},\cite{SBJ},\cite{FE}) to solve the equation $f(x)=0$ in the interval $[a,b]$, we find a relation among the suitable point $x_{n-1}$, the approximate root $x_{n}$ and the function derivative at $x=x_{n-1}$ so that $x_{n-1}$ and $x_{n}$ both belong to $[a,b]$ as follows:
\begin{eqnarray} 
x_{n}=x_{n-1}-\frac{f(x_{n-1})}{f'(x_{n-1})}
\end{eqnarray} 
But, the Newton's method takes a long step since reaching a solution as $x_{n}$ might create $n$ steps. In this section, we present a variant of Halley and Newton's method for reaching an approximation solution as $x^{*}\in[a,b]$ with at least possible steps provided that we firstly find a suitable point $x_{0}\in[a,b]$. This method is stated in the following theorem:

\begin{figure}[h!]
\centering
  \includegraphics[width=1\linewidth]{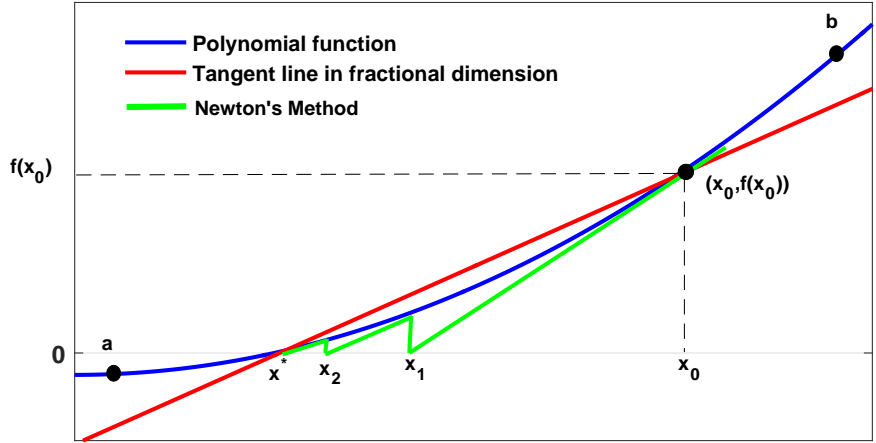}
  \caption{Fractional and classical derivatives at the point $x_{0}$ and solution $x^{*}$ by fractional derivative in comparison with Newton's method} \label{Fig:1}
\end{figure}
A variant of Halley and Newton's method along with fractional derivative is stated by theorem1.\\
\textbf{Theorem 1}\\
Let the function $f(x)$ be continuous and differentiable up to $nth$-order in the interval $[a,b]$ and according to Fig.1, if the root $x^{*}\in[a,b]$ denotes a solution for $f(x)=0$ and $x_{0}\in[a,b]$ an arbitrarily chosen  point, then 
\begin{eqnarray} 
x^{*}=x_{0}-\frac{f(x_{0})}{D^{\alpha}f(x)|_{x=x_{0}}}
\end{eqnarray} 
and
\begin{eqnarray} 
|\frac{f'(x_{0})\pm \sqrt{f'(x_{0})^{2}-2f(x_{0})f^{''}(x_{0})}}{f^{''}(x_{0})}|<1
\end{eqnarray}
provided that $|\frac{f(x_{0})}{D^{\alpha}f(x)|_{x=x_{0}}}|<1$,
where $D^{\alpha}f(x)|_{x=x_{0}}$ denotes the fractional derivative of $f(x)$ of the order $\alpha$ at the point $x_{0}$.\\
\textbf{\textit{Proof}}\\
Since the function $f(x)$ is continuous and differentiable up to $nth$-order in the interval $[a,b]$, we can consider the function is identical to the Taylor's expansion for every point of this interval. Making Taylor's expansion at the point $x_{0}$, we have
\begin{eqnarray} 
f(x)=f(x_{0})+\frac{f'(x_{0})}{1!}(x-x_{0})+\frac{f''(x_{0})}{2!}(x-x_{0})^{2}+\dots
\end{eqnarray} 
Replacing $x=x^{*}$ in (2.4), yields
\begin{eqnarray} 
f(x^{*})=f(x_{0})+\frac{f'(x_{0})}{1!}(x^{*}-x_{0})+\frac{f''(x_{0})}{2!}(x^{*}-x_{0})^{2}+\dots=0
\end{eqnarray} 
According to Fig.1
\begin{eqnarray} 
\frac{f(x_{0})-0}{|x_{0}-x^{*}|}=D^{\alpha}f(x)|_{x=x_{0}}
\end{eqnarray}
where denotes the slope of the crossing line between $x_{0}$ and $ x^{*}$ created by fractional derivative at the point $x_{0}$.\\
Substituting $x^{*}-x_{0}=\frac{-f(x_{0}}{D^{\alpha}f(x)|_{x=x_{0}}}$ from (2.6) for (2.5) and considering the $|\frac{-f(x_{0}}{D^{\alpha}f(x)|_{x=x_{0}}}|<1$ and truncating the fourth terms onward, we have
\begin{eqnarray} 
f(x_{0})-f'(x_{0})\left\{\frac{-f(x_{0})}{D^{\alpha}f(x)|_{x=x_{0}}}\right\}+\frac{f''(x_{0})}{2}\left\{\frac{-f(x_{0})}{D^{\alpha}f(x)|_{x=x_{0}}}\right\}^{2}\approx 0
\end{eqnarray} 
Letting $Z=\frac{-f(x_{0})}{D^{\alpha}f(x)|_{x=x_{0}}}$ in (2.7), we find a quadratic equation
\begin{eqnarray} 
f''(x_{0})Z^{2}-2f'(x_{0})Z+2f(x_{0})\approx 0
\end{eqnarray}
Solving this equation with regard to $Z$ gives us the inequality (2.3) and the solution $x^{*}$ given by (2.2). \\ Some examples of the quadratic and cubic polynomials are stated for representing the accuracy and speed of the mentioned method in theorem1.
\textbf{Example 1}\\
Solve the quadratic equation $x^{2}+3x+1=0$.\\
Let the function $f(x)=x^{2}+3x+1$, then $f'(x)=2x+3$ and $f''(x)=2$.\\
Letting arbitrarily $x_{0}=5$ in $f'$, $f''$ and (2.3), we have\\
\begin{eqnarray} 
|\frac{13\pm \sqrt{13^{2}-2\times2\times41}}{2}|=\frac{13\pm2.36}{2}>1
\end{eqnarray}
Thus, $x_{0}=5$ is not a suitable one. 
But, if we let $x_{0}=0.5$, then $f(0.5)=2.75$, $f'(0.5)=4$ and $f''(0.5)=2$, then
\begin{eqnarray} 
|\frac{4\pm \sqrt{4^{2}-2\times2.75\times2}}{2}|=\frac{4\pm2.24}{2}
\end{eqnarray}
As we see, only one of the solutions can be correct. The corrected solution is: $\frac{4-2.24}{2}=0.881<1$, then using (2.2) and (2.3), we find $x^{*}=0.5-0.881=-0.381$, therefore, one of the solutions is: $x^{*}=-0.381$.\\
If we let $x_{0}=-3.5$, then $f(-3.5)=2.75$, $f'(-3.5)=-4$ and $f''(-3.5)=2$, then we will find the other solution as follows:\\
$x^{*}=-3.5-(-0.881)=-2.619$. Thus, the other solution is: $x^{*}=-2.619$.
\textbf{Example 2}\\  
Solve the cubic equation $x^{3}+2x^{2}-4x-8=0$.\\ 
Let the function $f(x)=x^{3}+2x^{2}-4x-8$ and $f'(x)=3x^{2}+4x-4$, $f''(x)=6x+4$.\\
Let $x_{0}=3$, then $f(3)=25$, $f'(3)=35$, and $f''(3)=22$. Regarding (2.3), we have 
\begin{eqnarray} 
|\frac{35\pm \sqrt{35^{2}-2\times25\times22}}{22}|
\end{eqnarray}
where is not acceptable since under the square root symbol is a negative number.\\
Choosing $x_{0}=2.5$ and substituting for the function $f(x)$, we have\\
$f(2.5)=10.125$, $f'(2.5)=24.75$ and $f''(2.5)=19$. (2.3) gives
\begin{eqnarray} 
|\frac{24.75\pm \sqrt{24.75^{2}-2\times10.125\times19}}{19}|
\end{eqnarray} 
The only the solution $|\frac{24.75- \sqrt{24.75^{2}-2\times10.125\times19}}{19}|=0.508<1$ is acceptable. Therefore, regarding (2.2) and (2.3), the solution for cubic function $f(x)$ is:\\
$x^{*}=x_{0}-0.508=2.5-0.508\approx 2$.\\
The other solution is obtained by putting $x_{0}=-2$, where gives:\\
$f(-2)=0$, $f'(-2)=0$ and $f''(-2)=-8$ and from (2.3) we have:\\
$x^{*}=x_{0}=-2$.\\
Just, we are ready to calculate the value of fractional order given in theorem 1 as theorem 2.\\\\
\textbf{Theorem 2}\\
Let $_{RL}D^{\alpha}f(x)|_{x=x_{0}}$ be Riemann-Liouville fractional derivative of the function of the order $\alpha$ applied for (2.3) and theorem1's conditions hold, then value of $\beta=-\alpha$ is obtained:
\begin{eqnarray} 
\beta=\frac{\log |\Gamma(\beta+1)|+\log  \left\{ |\frac{f(x_{0})f^{''}(x_{0})}{f'(x_{0})\pm \sqrt{f'(x_{0})^{2}-2f(x_{0})f^{''}(x_{0})}}|\right\}-\log |f(\xi_{1})|}{\log(x_{0}-a)}
\end{eqnarray}  
where $a<\xi_{1}<x_{0}$.\\
\textbf{\textit{Proof}}\\
If we make use of the Riemann-Liouville derivative $_{RL}D^{\alpha}f(x)|_{x=x_{0}}$, we have
\begin{eqnarray} 
_{RL}D^{\alpha}f(x)|_{x=x_{0}}=\frac{1}{\Gamma(n-\alpha)}\frac{d^{n}}{dx^{n}}\int_{a}^{x}(x-\xi)^{n-\alpha-1}f(\xi)d\xi|_{x=x_{0}}
\end{eqnarray}  
where $n-1\leq \alpha<n$.\\
Leibniz integral rule states
\begin{eqnarray} 
\frac{d}{dx}\int_{a(x)}^{b(x)}f(x,t)dt=\int_{a(x)}^{b(x)}\frac{\partial}{\partial x}f(x,t)dt+f(x,b(x))\frac{d b(x)}{dx}-\nonumber\\f(x,a(x))\frac{d a(x)}{dx}~~~~~~~~~~~~~~~~~~~~~~~~~~~
\end{eqnarray}   
Let $f(x,t)=(x-\xi)^{n-\alpha-1}f(\xi)$, $a(x)=a$, and $b(x)=x$ in (2.15), and taking $n$ times derivation of (2.15), we have (2.14) as follows:
\begin{eqnarray} 
_{RL}D^{\alpha}f(x)|_{x=x_{0}}=\frac{1}{\Gamma(-\alpha)}\int_{a}^{x}(x-\xi)^{-\alpha-1}f(\xi)d\xi|_{x=x_{0}}
\end{eqnarray} 
Let $\beta=-\alpha$ at the point $x=x_{0}$ and substituting for (2.3), then\\
\begin{eqnarray} 
_{RL}D^{-\beta}f(x)|_{x=x_{0}}=\frac{1}{\Gamma(\beta)}\int_{a}^{x_{0}}(x_{0}-\xi)^{\beta-1}f(\xi)d\xi=\nonumber\\ \frac{f(x_{0})f^{''}(x_{0})}{f'(x_{0})\pm \sqrt{f'(x_{0})^{2}-2f(x_{0})f^{''}(x_{0})}}~~~~~~~~~~~~
\end{eqnarray} 
Regarding mean value theorem in fractional integrals for (2.17), yields
\begin{eqnarray} 
\frac{1}{\Gamma(\beta)}\int_{a}^{x_{0}}(x_{0}-\xi)^{\beta-1}f(\xi)d\xi=\frac{f(\xi_{1})}{\Gamma(\beta)}\int_{a}^{x_{0}}(x_{0}-\xi)^{\beta-1}d\xi=\nonumber\\ \frac{f(\xi_{1})}{\Gamma(\beta)}\frac{(x_{0}-a)^{\beta}}{\beta}=\frac{f(\xi_{1})}{\Gamma(\beta+1)}(x_{0}-a)^{\beta}~~~~~~~~~~
\end{eqnarray} 
where $a<\xi_{1}<x_{0}$.\\
Thus,\\
Just, (2.17) and (2.18) imply that
\begin{eqnarray} 
\frac{f(\xi_{1})}{\Gamma(\beta+1)}(x_{0}-a)^{\beta}=\frac{f(x_{0})f^{''}(x_{0})}{f'(x_{0})\pm \sqrt{f'(x_{0})^{2}-2f(x_{0})f^{''}(x_{0})}}
\end{eqnarray} 
manipulating (2.19) for obtaining $\beta$ and taking logarithm from both left and right sides of (2.19) gives our desirable result as (2.13).\\
Examples 3 and 4 are presented here to make clear the method of obtaining $\alpha$ or $\beta$ mentioned in theorem 2.\\ 
\textbf{Example 3}\\ 
Referring example 1 and considering the solution $x_{0}=0.5$ and $a=0$, we arbitrarily choose a number between 0 and 0.5 as $\xi_{1}=0.25$, then $f(0.25)=1.812$ and $\frac{f(x_{0})f^{''}(x_{0})}{f'(x_{0})\pm \sqrt{f'(x_{0})^{2}-2f(x_{0})f^{''}(x_{0})}}=3.12$. $\beta$ can be calculated by (2.13) as follows:
\begin{eqnarray} 
\beta=-1.443\log \Gamma(\beta+1)-0.784
\end{eqnarray}  
Just, we solve the equation (2.20) via a maple16 computer program of approximation less than 0.013 and the initial $\beta=-2.01$ as follows:\\
Restart:\\ 
Bet := -2.01;\\ 
for i to 99 do\\ 
J := Bet+1;\\ 
N := evalf(GAMMA(J));\\ 
M := -1.443*ln(abs(N))-.784;\\ 
if (abs(Bet-M) $\leq$ 0.013) then print(Bet);\\ print("**************************************");\\ 
else Bet := -2.01-(1/100)*i;\\ 
end if;\\ 
end do;\\
The final solution is: $\beta=-2.86$\\
Note that choosing $\beta=-2.01$ is due to looking at (2.20), we intuitively understand that $\beta+1$ puts between -2 and -1. In such a case, $\log(\Gamma(\beta+1))$ gets a positive value since $\Gamma(\beta+1)\geq 1$ and (2.20) might have a solution.\\
For the other solution in example 1,  
let $x_{0}=-3.5$,  choose $a=-5$ and $\xi_{1}=\frac{-3.5+(-5)}{2}=-4.25$.\\ Let $f(-3.5)=2.75$, $f'(-3.5)=-4$, $f''(-3.5)=2$, $f(-4.25)=6.3125$, then (2.13) is obtained as follows:
\begin{eqnarray} 
\beta=2.466\log \Gamma(\beta+1)-1.738
\end{eqnarray}    
Choosing the initial value for $\beta=-4.01$, the final solution of the approximation 0.013  is: $\beta=-4.38$. For calculating the solution, we make use of Maple16 as follows:\\
Restart:\\ 
Bet := -4.01;\\ 
for i to 99 do\\ 
J := Bet+1;\\ 
N := evalf(GAMMA(J));\\ 
M := 2.466*ln(abs(N))-1.738;\\ 
if (abs(Bet-M) $\leq$ 0.013) then print(Bet);\\ print("**************************************");\\ 
else Bet := -4.01-(1/100)*i\\ 
end if;\\ 
end do;\\
\textbf{Example 4}\\ 
Referring example 2 and choosing $x_{0}=2.5$ as an initially right point for solution $x^{*}$, then choosing $a=2$ and $\xi_{1}=\frac{2.5+2}{2}=2.25$, we find $f(2.5)=10.125$, $f'(2.5)=24.75$ and $f''(2.5)=19$, $f(2.25)=4.516$ and $\beta$ is obtained by the following relation: 
\begin{eqnarray} 
\beta=-1.442\log \Gamma(\beta+1)-2.1418
\end{eqnarray} 
Choosing the initial value for $\beta=-3.01$ of the approximation less than or equal to 0.009, we find the final solution $\beta=-3.21$ via the following Maple 16 program:\\
Restart:\\ 
Bet := -3.01;\\
for i to 99 do\\
J := Bet+1;\\ 
N := evalf(GAMMA(J));\\ 
M := -1.442*ln(abs(N))-2.1418;\\ 
if (abs(Bet-M) $\leq$ 0.009) then print(Bet);\\
print("**************************************");\\
else Bet := -3.01-(1/100)*i\\
end if;\\
end do;\\
Regarding the other root $x_{0}=-2$, we could not find any $\beta$.
\subsection{\textbf{The functions of several variables-A generalization to \textbf{$\mathbb{R}^{n}$}}}
Let $X^{*}=(x^{*}_{1},\dots,x^{*}_{n})$ be a root of the convex function $F(X)=F(x_{1},\dots,x_{n})$ and the initial point at $X_{0}=(x_{01},\dots,x_{0n})$, then $F(X^{*})=0$. A generalization of theorem 1 only for convex functions is stated as theorem 3.\\
\textbf{Theorem 3}\\
Given the convex function $F(x_{1},\dots,x_{n})$, if the root $X^{*}\in[(a_{1},...,a_{n}),(b_{1},...,b_{n})]$ denotes a solution for $F(X)=F(x_{1},\dots,x_{n})=0$ and $X_{0}=(x_{01},\dots,x_{0n})\in [(a_{1},...,a_{n}),(b_{1},...,b_{n})]$ an arbitrary point of the function, then 
\begin{eqnarray} 
\frac{F(X_{0})-0}{|X_{0}-X^{*}|}=(\nabla^{\alpha}.\vec{h})F(X)|_{X=X_{0}}
\end{eqnarray}
and 
\begin{eqnarray} 
F(X)=F(X_{0})+(h_{1}D_{1}+\dots+h_{n}D_{n})F(X_{0})+\nonumber\\ \frac{1}{2!}(h_{1}D_{1}+\dots+h_{n}D_{n})^{2}F(X_{0})+\dots+\nonumber\\ \frac{1}{(m-1)!}(h_{1}D_{1}+\dots+h_{n}D_{n})^{m-1}F(X_{0})+\dots
\end{eqnarray}
where  
$D_{1}=\frac{\partial}{\partial x_{1}},\dots,D_{n}=\frac{\partial}{\partial x_{n}}$,     
$\vec{h}=(h_{1},\dots,h_{n})$ denotes the unit vector,  $\nabla^{\alpha}=\frac{\partial^{\alpha}}{\partial x_{1}}i_{1}+\frac{\partial^{\alpha}}{\partial x_{2}}i_{2}+\dots+\frac{\partial^{\alpha}}{\partial x_{n}}i_{n}$ and $i_{1},\dots,i_{n}$ denote the vector units of the space and $\alpha$ order of fractional partial derivative of gradient $\nabla^{\alpha}$\\
\textbf{\textit{Proof}}\\ 
By the Chorlton's paper \cite{CF}, we find (2.24) and by theorem 1, we find (2.23) since
\begin{eqnarray} 
F(X^{*})=F(X_{0})+(h^{*}_{1}D_{1}+\dots+h^{*}_{n}D_{n})F(X_{0})+\nonumber\\ \frac{1}{2!}(h^{*}_{1}D_{1}+\dots+h^{*}_{n}D_{n})^{2}F(X_{0})+\dots+\nonumber\\ \frac{1}{(m-1)!}(h^{*}_{1}D_{1}+\dots+h^{*}_{n}D_{n})^{m-1}F(X_{0})+\dots=0
\end{eqnarray}
and
\begin{eqnarray} 
h^{*}_{1}=(x^{*}_{1},\dots,x^{*}_{n})-(x_{01},x_{2},\dots,x_{n})
\end{eqnarray}
\begin{eqnarray} 
h^{*}_{2}=(x^{*}_{1},\dots,x^{*}_{n})-(x_{1},x_{02},\dots,x_{n})
\end{eqnarray}
~~~~~~~~~~~~~~~~~~~~~~~~~~~~~\vdots~~~~~~~~~~~~~~~~~\vdots~~~~~~~~~~~~~~~~~\vdots
\begin{eqnarray} 
h^{*}_{n}=(x^{*}_{1},\dots,x^{*}_{n})-(x_{1},x_{2},\dots,x_{0n})
\end{eqnarray}
\textbf{Theorem 4}\\
Let us functions of several variables be analogues to the relations (2.2), (2.3),  (2.6), (2.25), (2.7) and (2.8) and they hold for functions of several variables as given $F(X)$ in theorem3 provided that 
\begin{eqnarray} 
h^{*}_{1}=h^{*}_{2}=\dots=h^{*}_{n}
\end{eqnarray}
and 
\begin{eqnarray} 
dx_{1}=dx_{2}=\dots=dx_{n}
\end{eqnarray}   
then, we have
\begin{eqnarray} 
DF(X)|_{X_{0}}=\sum_{i=1}^{n}\frac{\partial F}{\partial x_{i}}|_{X_{0}}
\end{eqnarray}
and  
\begin{eqnarray} 
D^{2}F(X)|_{X_{0}}=\sum_{i=1}^{n}\frac{\partial^{2} F}{\partial x_{i}^{2}}|_{X_{0}}+2\sum_{i,j}\frac{\partial^{2} F}{\partial x_{i}\partial x_{j}}|_{X_{0}}
\end{eqnarray}
where $DF(X)|_{X_{0}}$ denotes the first derivative of the function of several variables $F(X)$ at the point $X_{0}$ and $D^{2}F(X)|_{X_{0}}$ denotes the second derivative of same function  $F(X)$ at same point.\\ 
Like (2.2) and (2.3), we have
\begin{eqnarray} 
X^{*}=X_{0}-\frac{F(X_{0})}{(\nabla^{\alpha}.\vec{h})F(X)|_{X=X_{0}}}
\end{eqnarray}
and
\begin{eqnarray} 
|\frac{DF(X)|_{X_{0}}\pm \sqrt{(DF(X)|_{X_{0}})^{2}-2F(X_{0})D^{2}F(X)|_{X_{0}}}}{D^{2}F(X)|_{X_{0}}}|<1
\end{eqnarray}
provided that $|\frac{F(X_{0})}{(\nabla^{\alpha}.\vec{h})F(X)|_{X=X_{0}}}|<1$.\\
\textbf{\textit{Proof}}\\ 
The proof is trivial since the relation (2.24) and condition (2.29) regarding chain rule says us
\begin{eqnarray} 
dF(X)=\sum_{i=1}^{n}\frac{\partial F}{\partial x_{i}}dx_{i}|_{X_{0}}
\end{eqnarray}
Dividing both hand-sides by $dx_{i}$ and considering the condition (2.30), we find (2.31). Similarly, the relation (2.32) is proved. (2.33) and (2.34) are exactly similar to (2.2) and (2.3) as well. The conditions (2.29),(2.30) and (2.34) are the necessary for solution of $F(X)=0$ not sufficient ones.
To learn more about convex cone and fractional dimension hyperplanes, refer to  the paper\cite{MA2}\\
\textbf{Example 5}\\ 
Solve the nonlinear system of equations 
\[  \left\{ \begin{array}{c}
3x_{1}-cos(x_{2}x_{3})-\frac{3}{2}=0~~~~~~~~~\\
4x_{1}^{2}-625x_{2}^{2}+2x_{2}-1=0~~~~~~\\
\exp(-x_{1}x_{2})+20x_{3}+\frac{10\pi-3}{3}=0
\end{array} \right.  \]
\begin{eqnarray} 
\end{eqnarray}
Let us introduce the associated function
\[ G(x)=\left(\begin{array}{cc}
3x_{1}-cos(x_{2}x_{3})-\frac{3}{2}\\
4x_{1}^{2}-625x_{2}^{2}+2x_{2}-1\\
\exp(-x_{1}x_{2})+20x_{3}+\frac{10\pi-3}{3}
\end{array} \right) \] 
\begin{eqnarray} 
\end{eqnarray}
where 
\[X=\left[\begin{array}{cc}
x_{1}\\
x_{2}\\
x_{3}
\end{array} \right] \]
\begin{eqnarray} 
\end{eqnarray} 
One might now define the objective function
\begin{eqnarray} 
F(X)=\frac{1}{2}G^{T}(X)G(X)=~~~~~~~~~~~~~~~\nonumber\\ \frac{1}{2}\{(3x_{1}-cos(x_{2}x_{3})-\frac{3}{2})^{2}+(4x_{1}^{2}-625x_{2}^{2}+2x_{2}-1)^{2}+\nonumber\\ (\exp(-x_{1}x_{2})+20x_{3}+\frac{10\pi-3}{3})^{2}\}~~~~~~~~~~~
\end{eqnarray}
which we will attempt to minimize.\\ 
Regarding (2.31), we have
\begin{eqnarray} 
DF(X)|_{X_{0}}=\frac{\partial F}{\partial x_{1}}|_{X_{0}}+\frac{\partial F}{\partial x_{2}}|_{X_{0}}+\frac{\partial F}{\partial x_{3}}|_{X_{0}}=~~~~~~~\nonumber\\
\{3+x_{3}sin(x_{2}x_{3})+x_{2}sin(x_{2}x_{3})\}\{3x_{1}-cos(x_{2}x_{3})-1.5\}|_{X_{0}}+\nonumber\\ \{8x_{1}-1250x_{2}+2\}\{4x_{1}^{2}-625x_{2}^{2}+2x_{2}-1\}|_{X_{0}}+~~~~~~~\nonumber\\ \{-(x_{1}+x_{2})\exp(-x_{1}x_{2})+20\}\{\exp(-x_{1}x_{2})+20x_{3}+\frac{10\pi-3}{3}\}|_{X_{0}}
\end{eqnarray}
and regarding (2.32)
\begin{eqnarray} 
D^{2}F(X)|_{X_{0}}=~~~~~~~~~~~~~~~~~~~~~~~~~~~~~~~~~~~~~~\nonumber\\ \frac{\partial^{2} F}{\partial x_{1}^{2}}|_{X_{0}}+\frac{\partial^{2} F}{\partial x_{2}^{2}}|_{X_{0}}+\frac{\partial^{2} F}{\partial x_{3}^{2}}|_{X_{0}}+2\frac{\partial^{2} F}{\partial x_{1}\partial x_{2}}|_{X_{0}}+2\frac{\partial^{2} F}{\partial x_{1}\partial x_{3}}|_{X_{0}}+2\frac{\partial^{2} F}{\partial x_{2}\partial x_{3}}|_{X_{0}}=\nonumber\\
\{x_{2}^{2}cos(x_{2}x_{3})+x_{3}^{2}cos(x_{2}x_{3})+2sin(x_{2}x_{3})+2x_{2}x_{3}cos(x_{2}x_{3})\}\times~~~~~~\nonumber\\ 
\{3x_{1}-cos(x_{2}x_{3})-\frac{3}{2}\}|_{X_{0}}+\{3+x_{2}sin(x_{2}x_{3})+x_{3}sin(x_{2}x_{3})\}^{2}|_{X_{0}}-~~~\nonumber\\
1242\{4x_{1}^{2}-625x_{2}^{2}+2x_{2}-1\}|_{X_{0}}+\{8x_{1}-1250x_{2}+2\}^{2}|_{X_{0}}+~~~~~~~~~~\nonumber\\ \{x_{1}^{2}\exp(-x_{1}x_{2})+x_{2}^{2}\exp(-x_{1}x_{2})-2\exp(-x_{1}x_{2})+2x_{1}x_{2}\exp(-x_{1}x_{2})\}\times\nonumber\\ \{\exp(-x_{1}x_{2})+20x_{3}+\frac{10\pi-3}{3}\}|_{X_{0}}+\{-(x_{1}+x_{2})\exp(-x_{1}x_{2})+20\}^{2}|_{X_{0}}~~~
\end{eqnarray}
Just we guess a initial point so that
\[X^{(0)}=\textbf{0}=\left[\begin{array}{cc}
0\\
0\\
0
\end{array} \right] \]
\begin{eqnarray} 
\end{eqnarray} 
but, the computations show us that $Delta=(DF(X)|_{X_{0}})^{2}-2F(X_{0})D^{2}F(X)|_{X_{0}}<0$ expressed in the relation (2.34) is negative and there is no solution. Therefore, we need to guess another initial point. To find better point, we analyze the system (2.36) and find the following point, which seems to be close to the solutions since the first equation says us that $\frac{1}{6}\leq x_{1}\leq \frac{5}{6}$:
\[X^{(0)}=\left[\begin{array}{cc}
\frac{1}{2}\\
0.0032\\
-0.523
\end{array} \right] \]
\begin{eqnarray} 
\end{eqnarray}    
again, we have $Delta<0$. The other guess we make is only a change in $x_{1}$ as follows:
\[X^{(0)}=\left[\begin{array}{cc}
\frac{2}{3}\\
0.0032\\
-0.523
\end{array} \right] \]
\begin{eqnarray} 
\end{eqnarray}   
and we get 
\begin{eqnarray}
F(\frac{2}{3},0.0032,-0.523)=0.427 
\end{eqnarray}
This states that we are closing to the minimum of the objective function $F(X)$ and finally closing to the solution.\\
Just, using the relation (2.33), we find two groups of the solutions as follow:
\[X^{(1)}=\left[\begin{array}{cc}
0.6302\\
-0.0332\\
-0.5594
\end{array} \right] \]
\begin{eqnarray} 
\end{eqnarray}   
and 
\[X^{(1)}=\left[\begin{array}{cc}
0.707\\
0.0441\\
-0.482
\end{array} \right] \]
\begin{eqnarray} 
\end{eqnarray}   
For the first group, we have
\begin{eqnarray}
F(0.6302,-0.0332,-0.5594)=0.44 
\end{eqnarray}
and for second one
\begin{eqnarray}
F(0.707,-0.0441,-0.482)=0.4 
\end{eqnarray}
If we change the first variable of the initial guess (2.44) by the first one of the first above group, we find
\[X^{(1)}=\left[\begin{array}{cc}
0.6302\\
0.0032\\
-0.523
\end{array} \right] \]
\begin{eqnarray} 
\end{eqnarray}
and 
\begin{eqnarray}
F(0.6302,0.0032,-0.523)=0.358 
\end{eqnarray}     
This means that we are closing much more to the solution.\\
\textbf{A comparison of our method with Gradient descent one}\\
If we get started with the initial guess given in (2.42) for Gradient Descent method, we find
\begin{eqnarray}
X^{(1)}=\textbf{0}-\gamma_{0}\Delta F(\textbf{0})=\textbf{0}-\gamma_{0}J_{G}(\textbf{0})G(\textbf{0}) 
\end{eqnarray}
where $J_{G}$ denotes Jacobian matrix and $G(\textbf{0})$ denotes matrix (2.37) at the vector point $\textbf{0}$. Jacobian matrix is given by

\[ J_{G}(X)=\left(\begin{array}{cc}
3~~~~~~~~~~x_{3}sin(x_{2}x_{3})~~~~~~~~~~x_{2}sin(x_{2}x_{3})\\
8x_{1}~~~~~~~~~~~~~~~~~2-1250x_{2}~~~~~~~~~~~~~~~~~0\\
-x_{2}\exp(-x_{1}x_{2}~~~-x_{1})\exp(-x_{1}x_{2})~~~20
\end{array} \right) \] 
\begin{eqnarray} 
\end{eqnarray}
then calculating at the point $\textbf{0}$, we have
\[ J_{G}(\textbf{0})=\left(\begin{array}{cc}
3~~~~~~~~~0~~~~~~~~~~0\\
0~~~~~~~~~2~~~~~~~~~~0\\
0~~~~~~~~~0~~~~~~~~~20
\end{array} \right) \] 
\begin{eqnarray} 
\end{eqnarray}
 and \[ G(\textbf{0})=\left(\begin{array}{cc}
-2.5\\
-1\\
10.472
\end{array} \right) \]
\begin{eqnarray} 
\end{eqnarray}
that implies that regarding (2.39), $F(\textbf{0})=58.456$. We see that the initial guess of the Gradient descent method is much more than to that of our method. Just for using (2.52), we need to find $\gamma_{0}$ as well as the initial guess (here we guess $\gamma_{0}=0.001$). Thus,
\[X^{(1)}=\left[\begin{array}{cc}
0.0075\\
0.002\\
-0.20944
\end{array} \right] \]
\begin{eqnarray} 
\end{eqnarray}   
where implies that $F(X^{(1)})=23.306$. The calculations show that we need to do 83 iterations of Gradient descent to reach a reasonably minimum value of the function $F(X)$. But, using our method , we could find a reasonable optimization at 5 iterations since our initial guess was much more close to the solution. In our method, not being negative Delta and condition (2.34) play a vital role in finding faster solutions.\\
If we get started with initial guess given by (2.44) using Gradient descent, then applying $\gamma_{0}=0.001$ we find
\[X^{(1)}=\left[\begin{array}{cc}
0.66402\\
0.00476\\
-0.52319
\end{array} \right] \]
\begin{eqnarray} 
\end{eqnarray}
and
\begin{eqnarray}
F(0.66402,0.00476,-0.52319)=0.417 
\end{eqnarray} 
that after two iterations, we reach 0.417. In comparison with (2.51) our method shows that the function $F(X)$ gets values less than the Gradient descent one after two iterations.\\   
\textbf{Discussion}\\
As we stated before, to find a suitable $x_{0}$ for solving an equation by theorem1, we should make use of intuitive feelings. For instance in example1, if we look at the equation $x^{2}+3x+1=0$, we see that a solution must be negative. Certainly, we must find a $x_{0}$ closer than to negative numbers not positive ones. $x_{0}=5$ is farther than a negative solution but $x_{0}=0.5$ is better and closer than to it. Finding a suitable $x_{0}$ before using (2.2) and (2.3) makes the procedure of solution faster. For instance in example 2, we should be able to guess an approximation solution and trade off looking at equation $x^{3}+2x^{2}-4x-8=0$. Trading off we understand that $x_{0}=2.5$ is closer to actual solution $x^{*}=2$.\\
In examples 3 and 4, before using (2.13) and taking an initial $\beta$, one should pay attention to the signs of $f(\xi_{1})$, $\Gamma(\beta+1)$, and $\frac{f(x_{0})f^{''}(x_{0})}{f'(x_{0})\pm \sqrt{f'(x_{0})^{2}-2f(x_{0})f^{''}(x_{0})}}$ in the relation (2.19). Because, choosing $\beta$ affects to the sign of $\Gamma(\beta+1)$. For instance, if $-1<\beta<0$,  $\Gamma(\beta+1)>0$ or if  $-2<\beta<-1$, then $\Gamma(\beta+1)<0$ and so on.\\
In example 4, to obtain a $\beta$ for the solution $x_{0}=-2$, we could not find such $\beta$ since the term
\begin{eqnarray} 
\frac{f(x_{0})f^{''}(x_{0})}{f'(x_{0})\pm \sqrt{f'(x_{0})^{2}-2f(x_{0})f^{''}(x_{0})}}=\frac{0}{0}
\end{eqnarray} 
when $x_{0}=-2$ becomes indeterminate. After removing indeterminable term by evaluating its limit when $x$ tends to -2, we find the value 0, but there is no $\beta$ for this value 0.\\
As we understand from the paper,for having solutions close to real one, roots should be near to real solutions  either in the Gradient descent method or in our method. On the other hand, in Gradient descent method, we need find a coefficient $\gamma_{0}$ as well as initially guessed solutions. In our method, there is a criterion for filtering initially guessed solutions so that the number of iterations is limited with regard to Gradient descent method. Both methods show that initial guess should be close to real solutions to reach convergence. If the function is of several variables, then to obtain the best solution of least time, we should try to fit to inequality (2.34) along with closest guesses to real solutions.   
\section{Conclusions}
In this paper, we can replace the Newton's method applied to functions of one variable to obtain roots of an equation by the method stated in theorem1, which is so easier and faster than Newton's. We also give the method for finding order of a Riemann-Liouville fractional derivative so that the solution of an equation can be done so easier. The methods given by theorems 1 and 2 might be used for optimization as well.\\
The method given here can be extended and generalized to functions of several variables. Since theorems3 and 4  are the extensions of the theorem1 and 2, we can apply theorem1's method for functions of several variables for special cases.\\ 
\textbf{Acknowledgment}\\
We would like to thank the Editors and anonymous referees to review this paper and nice comments.


\newpage     
\begin{thebibliography}{99}
\bibitem{ACT} A. Akg\"{u}l, A. Cordero, J.R. Torregrosa, A fractional Newton method with
$2\alpha$th-order of convergence Appl. Math. Lett., 98 (2019) 344-351.

\bibitem{BA} A. Barvinok. A Course in Convexity, volume 54 of Graduate Studies in
Mathematics. Amer. Math. Soc., 2002.

\bibitem{CCT} G. Candelario, A. Cordero,  J.R.Torregrosa, Multipoint fractional iterative methods with ($2\alpha+1$)th-order of convergence for solving nonlinear problems, Math., 8 (2020) 1-15.

\bibitem{CF} F. Chorlton, Taylor's theorem for a function of several variables, Int.J.Math.Educ. Sci. Tech., 18, 2 (1987) 315-323.  

\bibitem{DBP} B.P. Demidovich, I.A. Maron, Computational mathematics, Medtech, U.K.,2017.

\bibitem{DFB} D.F. Bailey,A Historical Survey of Solution by Functional
Iteration,Math.Mag. 62, 3 (1989) 155-166. 

\bibitem{FE} G.E. Forsythe, Modern mathematics for the engineer,ed. E.F.Beckenbach, first ed., Chapt.17, what are relaxation methods?, 1956. 

\bibitem{GO} A.O. Gelfon, calculus of finite differences, Mir publications, Moscow, 1952.

\bibitem{HB} Harry Bateman, Halley's methods for solving equations, Amer. Math. Monthly 45 (1938)11-17.

\bibitem{MA2} A. Dorostkar, Relation between Roots and Tangent Lines of Function in Fractional Dimensions: A Method for Optimization Problems, Int. J. Math. Comput. Sci.13,9 (2019) 179-183.  

\bibitem{HBJ} J.B. Hiriart-Urruty and C. Lemare'chal, Convex Analysis and Minimization Algorithms, Springer, 1993. Two volumes.

\bibitem{HBJC}J.B. Hiriart-Urruty and C. Lemare'chal, Fundamentals of Convex Analysis.Springer, 2001. Abridged version of Convex Analysis and Minimization
Algorithms volumes 1 and 2.

\bibitem{LRS} S. R. Lay ,Convex Sets and Their Applications, John Wiley and Sons, 1982.

\bibitem{KV} V. Klee. ,What is a convex set?, Amer.Math.Month.,
78,6 (1971) 616-631.

\bibitem{SBJ} J.B. Scarborough, Numerical mathematical analysis, Chap XVIII, 1955.


\bibitem{SG} M.G. Salvadori, M.L. Baron, Numerical methods in engineering, Chap. 1, 1952.

\end{thebibliography}
\end{document}